\documentclass{birkjour}
\usepackage{amsthm}
\usepackage{amssymb}
\usepackage{amsmath}
\usepackage{mathtools}
\usepackage{verbatim}

\DeclarePairedDelimiter\floor{\lfloor}{\rfloor}

\newtheorem{theorem}{Theorem}[section]
\newtheorem{lemma}[theorem]{Lemma}
\newtheorem{kol}[theorem]{Corollary}
\theoremstyle{remark}
\newtheorem{bem}[theorem]{Remark} 
\newtheorem*{theorem*}{Theorem}

\def\bR{{\mathbb R}}
\def\bC{{\mathbb C}}
\def\te{{\textnormal e}}
\def\td{{\textnormal d}}
\def\ti{{\textnormal i}}
\def\opt{{\textnormal {opt}}}

\numberwithin{equation}{section}

\allowdisplaybreaks

\begin{document}

\title[A quantified Tauberian theorem for Laplace-Stieltjes transform]
 {A quantified Tauberian theorem for\\
 Laplace-Stieltjes transform}

\author{Markus Hartlapp}

\address{
Institute of Analysis\\
Technische Universit\"{a}t Dresden\\
Zellescher Weg 12-14\\
01069 Dresden\\
Germany}

\email{markus\_sebastian.hartlapp@tu-dresden.de}

\thanks{The author is supported by a doctoral scholarship of Konrad-Adenauer-Stifung} 

\subjclass{Primary 40E05; Secondary 40A05, 40E10, 44A10}

\keywords{Tauberian theorem, rates of decay, Laplace-Stieltjes transform}

\begin{abstract}
We prove a quantified Tauberian theorem involving Laplace-Stieltjes transform which is motivated by the work of Ingham and Karamata. 
For this, we consider functions which are locally of bounded variation and, therefore, get a generalisation of some results of Batty and Duyckaerts.
We show that our theorem can be applied to special Dirichlet series. 
\end{abstract}

\maketitle

\section{Introduction}

Considering Tauberian theorems which involve Laplace-Stieltjes transform is a business with a history of over 100 years. In 1916 Riesz generalised 
his observations about Dirichlet series to Laplace-Stieltjes transforms of functions which are locally of bounded variation \cite{Riesz16}. Some years 
later his work was refined independently by Ingham \cite{Ingham35} and Karamata \cite{Karamata34} who stated the following Tauberian theorem.

\begin{theorem}\label{ING}
Let $X$ be a Banach space, $A: [0,\infty) \to X$ locally of bounded variation, $A(0)=0$ and assume that there are $C'>0$ and $x_0 > 0$ so that

\begin{equation*}
 \limsup\limits_{t \to \infty} \left| \left| \textnormal{e}^{-x_0 t} \int_0^{t} \textnormal{e}^{x_0 s} \textnormal{d}A(s) \right| \right| \leq C'.
\end{equation*}
Then $f(z) = \int_0^{\infty} \textnormal{e}^{-zs} \textnormal{d}A(s)$ is convergent for every $z \in \mathbb{C}$ with $\rm{Re}(z) > 0$. 
Suppose further, that for some $A_{\infty} \in X$ the function $z \mapsto \frac{f(z)-A_{\infty}}{z}$ admits a continuous extension to the closed half-plane 
$\{z \in \bC \; | \; \textup{Re} (z) \geq 0\}$. Then

\[
 \limsup\limits_{t \to \infty} || A(t) - A_{\infty} ||  \leq 2C'.
\]
\end{theorem}
\noindent Actually, Ingham and Karamata showed this only for scalar-valued $A$ but there is no difficulty to adapt the proof for the vector-valued case.

In his approach to the Prime Number Theorem, Newman gave a new proof of Riesz' Tauberian theorem for Dirichlet series \cite{Newman80}. This proof was adapted to
Laplace transforms of bounded measurable functions by Korevaar \cite{Korevaar82} and Zagier \cite{Zagier97} who obtained special cases of the Ingham-Karamata-theorem. 
For an overview of the development of Tauberian theory see \cite{Korevaar04}.

The Newman-Korevaar-Zagier technique helped to state new kinds of results in the theory of stability of operator semigroups, see for example \cite{Arendt01}. 
Recent results in stability theory gave not only conditions for stability but stated convergence rates for semigroups $(T(t))_{t\geq0}$ and, analogously, for 
bounded measurable functions $f: [0,\infty) \to X$, where $X$ is a Banach space, for large times, see \cite{Batty08},\cite{Chill16}. 

In this paper we combine the ideas of Ingham and Karamata with those of Batty and Duyckaerts. Therefore, we get both a quantitative version of Theorem \ref{ING}
and a generalization of \cite[Theorem 4.1.]{Batty08} in the case $k=1$. In fact, we show the following.

\begin{theorem}\label{TH}
Take the same assumptions as in Theorem \ref{ING}. In addition, let $M: [0,\infty) \to [1,\infty)$ be a continuous, increasing function and 
$R: [0,\infty) \to [1,\infty]$ an increasing function. Assume that there exist $C>0, T \geq 0$ so that 

\begin{equation}\label{TC}
 \sup\limits_{t > T} \; \sup\limits_{x_0\leq x\leq R(t)} \left|\left| x \textnormal{e}^{-xt} \int_0^{t} \textnormal{e}^{xs} \textnormal{d}A(s) \right|\right| \leq C.
\end{equation}
Suppose further that

\begin{equation}\label{region}
\begin{split}
& \text{$f$ has an analytic extension into the region} \\
& Q \coloneqq \left\{ z \in \mathbb{C} \; \big| \; 0 \geq x > - \frac{1}{M(|y|)} \right\} 
\end{split}
\end{equation}
and that

\begin{equation}\label{growth}
\text{$||f(z)|| \leq M(|y|) $ holds throughout $Q$, where $z=x+\ti y$.}
\end{equation}
Then there exist $K>0, T'\geq 0$ so that

\begin{equation*}
 || A(t) - f(0) ||  \leq K \max \left\{ \frac{1}{M^{-1}_{\log} \left( \frac{t}{4} \right)}, \frac{1}{R(t)} \right\}
\end{equation*}
for every $t>T'$, where $M^{-1}_{\log}$ is the inverse of the function 
$M_{\log}$ defined by $M_{\log}(a) = M(a)( \log a + \log M(a) - \frac{1}{2} \log(5C))$ for $a \geq 1$. 

\end{theorem}

\noindent If $R(t)$, for increasing $t$, is growing quickly enough then we get exactly the same rate as in \cite{Batty08}, \cite{Chill16} -- see Remark \ref{bounded functions} below -- 
but for a wider class of functions, namely functions which are locally of bounded variation. Regarding the assumptions we remark the following.

\begin{enumerate}
 \item[(i)] In addition to Ingham and Karamata, we assume the Tauberian condition \eqref{TC}. There is a function $A$ such that this 6condtion is not true for $T=0$; 
  see Remark \ref{bem}.  
 \item[(ii)] The continuation property \eqref{region} and the growth condition \eqref{growth} are as in \cite{Batty08}. They ensure that we get a quantitative result. 
\end{enumerate}

\noindent In the following section we give three useful lemmas for the proof of Theorem \ref{TH} which we present in Section 3. Subsequently, we show that Theorem 
\ref{TH} includes the result from \cite{Batty08} and can be applied to Dirichlet series $f(z) = \sum_{n=1}^{\infty} b_n n^{-z}$ with a bounded sequence of coefficients 
$(b_n)_{n \in \mathbb{N}}$. For the rest of the article, we define $H \coloneqq \{ z \in \bC \; | \; \textup{Re}(z) > 0 \}$ and $\bR_+ \coloneqq [0,\infty)$.

\section{Preliminaries}

In this section we prove three lemmas which will be applied in the proof of Theorem \ref{TH}. In all of them we deduce different conclusions from the same condition,
namely a condition similar to \eqref{TC}. 

\begin{lemma}\label{lemma1}
Let $X$ be a Banach space, $A: \mathbb{R_+} \to X$ locally of bounded variation, $A(0)=0$ and suppose that

\begin{equation*}
 \sup\limits_{t \in \mathbb{R_+}} \left| \left| \textnormal{e}^{-xt} \int_0^{t} \textnormal{e}^{xs} \textnormal{d}A(s) \right| \right| \leq C
\end{equation*}
for some $x>0$ and some $C>0$. Let $z = x + \textnormal{i}y, y \in \bR $. Then 

\begin{equation*}
 \sup\limits_{t \in \mathbb{R_+}} \left| \left| \textnormal{e}^{-xt} \int_0^t \textnormal{e}^{zs} \textnormal{d}A(s) \right| \right| \leq C \left(1 + \frac{|y|}{x} \right).
\end{equation*}
\end{lemma}

\begin{proof}
First, define $G: \bR_+ \to X$ by $G(s) = \int_0^s \textnormal{e}^{xr} \textnormal{d}A(r), s \in \bR_+$. Then we have by properties of the Riemann-Stieltjes integral

\begin{equation*}
\left| \left| \textnormal{e}^{-xt} \int_0^t \textnormal{e}^{zs} \textnormal{d}A(s) \right| \right|
 = \left| \left| \textnormal{e}^{-xt} \int_0^t \textnormal{e}^{xs} \te^{\ti ys} \textnormal{d}A(s) \right| \right| 
 = \left| \left| \textnormal{e}^{-xt} \int_0^t \te^{\ti ys} \textnormal{d}G(s) \right| \right|. 
\end{equation*}
Integration by parts (cf. \cite[p.63]{Hille57}) and suitable estimates yield

\begin{align*}
\left| \left| \textnormal{e}^{-xt} \int_0^t \te^{\ti ys} \textnormal{d}G(s) \right| \right|
 & = \left| \left| \textnormal{e}^{-xt} \left( \left[ \te^{\ti ys} G(s) \right]_0^t - \ti y \int_0^t \te^{\ti ys} G(s) \td s \right) \right| \right| \\
 & \leq \left| \left| \textnormal{e}^{-xt} \te^{\ti yt} \int_0^t \te^{xr}\textnormal{d}A(r) \right| \right| \\
 & \quad \; + |y| \left| \left| \textnormal{e}^{-xt} \int_0^t \te^{\ti ys} \te^{xs} \te^{-xs} \left( \int_0^s \textnormal{e}^{xr} \textnormal{d}A(r) \right) \td s \right| \right|  \\
 & \leq C + |y| \te^{-xt} \int_0^t | \te^{\ti ys} | \te^{xs} \left|\left| \te^{-xs} \int_0^s \textnormal{e}^{xr} \textnormal{d}A(r) \right|\right| \td s \\
 & \leq C + C |y| \te^{-xt} \int_0^t \te^{xs} \td s \\
 & = C \left( 1 + \frac{|y|}{x} (1 - \te^{-xt}) \right) \leq C \left(1 + \frac{|y|}{x} \right).
\end{align*}
As this is true for every $t \in \bR_+$, we proved the claim.
\end{proof}

\begin{lemma}\label{lemma2}
Let $X$ be a Banach space, $A: \mathbb{R_+} \to X$ locally of bounded variation, $A(0)=0$ and suppose that

\begin{equation*}
 \sup\limits_{t \in \mathbb{R_+}} \left| \left| \textnormal{e}^{-xt} \int_0^{t} \textnormal{e}^{xs} \textnormal{d}A(s) \right| \right| \leq C
\end{equation*}
for some $x>0$ and some $C>0$. Let $z = x + \textnormal{i}y, y \in \bR $. Then 

\begin{equation}\label{l2}
 \sup\limits_{t \in \mathbb{R_+}} \left| \left| \textnormal{e}^{xt} \int_t^{\infty} \textnormal{e}^{-zs} \textnormal{d}A(s) \right| \right| \leq C \left(3 + \frac{|y|}{x} \right).
\end{equation}
\end{lemma}

\begin{proof}
Again we consider the function $G$ given by $G(s) = \int_0^s \textnormal{e}^{xr} \textnormal{d}A(r)$, $s \in \bR_+$. First, we show that the integral in \eqref{l2}
exists. We have

\begin{align*}
\left| \left| \int_t^v \textnormal{e}^{-zs} \textnormal{d}A(s) \right| \right| & = \left| \left| \int_t^v \te^{-\ti ys} \textnormal{e}^{-xs}  \textnormal{e}^{-xs} \textnormal{e}^{xs} \textnormal{d}A(s) \right| \right| \\
 & = \left| \left| \int_t^v \te^{-\ti ys} \textnormal{e}^{-2xs} \textnormal{d}G(s) \right| \right|.
\end{align*}
Integration by parts and suitable estimates yield

\begin{align*}
\left| \left| \int_t^v \te^{-\ti ys} \textnormal{e}^{-2xs} \textnormal{d}G(s) \right| \right|
 & = \left| \left| \left[ \te^{-\ti ys} \textnormal{e}^{-2xs} G(s) \right]_t^v + (2x +\ti y) \int_t^v \te^{-\ti ys} \textnormal{e}^{-2xs} G(s) \textnormal{d}s \right| \right| \\ 
 & \leq \left| \left| \te^{-\ti yv} \textnormal{e}^{-2xv} \int_0^v \textnormal{e}^{xr} \textnormal{d}A(r)  - \te^{-\ti yt} \textnormal{e}^{-2xt} \int_0^t \textnormal{e}^{xr} \textnormal{d}A(r)  \right| \right| \\
 & \quad \; + 2x \int_t^v \te^{-xs} \left| \left| \textnormal{e}^{-xs} \int_0^s \textnormal{e}^{xr} \textnormal{d}A(r)) \right| \right| \textnormal{d}s \\
 & \quad \; + |y| \int_t^v \te^{-xs} \left| \left| \textnormal{e}^{-xs} \int_0^s \textnormal{e}^{xr} \textnormal{d}A(r) \right| \right| \textnormal{d}s \\ 
 & \leq C \textnormal{e}^{-xv} + C \textnormal{e}^{-xt} + 2x C \int_t^v \te^{-xs} \td s + C |y| \int_t^v \te^{-xs} \td s \\
 & = C (\te^{-xv} + \te^{-xt}) - 2C(\te^{-xv} - \te^{-xt}) \\
 & \quad \; - C \frac{|y|}{x} (\te^{-xv} - \te^{-xt})
\end{align*}
which converges to $0$ for $v,t \to \infty$ and, therefore, the improper integral exists. For the estimate \eqref{l2} we write 

\begin{equation*}
\left| \left| \textnormal{e}^{xt} \int_t^{\infty} \textnormal{e}^{-zs} \textnormal{d}A(s) \right| \right| 
 = \lim\limits_{v \to \infty} \left| \left| \textnormal{e}^{xt}  \int_t^v \te^{-\ti ys} \textnormal{e}^{-2xs} \textnormal{d}G(s) \right| \right|.
\end{equation*}
The above estimate gives, for every $t \in \bR_+$,

\begin{align*}
\lim\limits_{v \to \infty} \left| \left| \textnormal{e}^{xt}  \int_t^v \te^{-\ti ys} \textnormal{e}^{-2xs} \textnormal{d}G(s) \right| \right|
 & \leq \lim\limits_{v \to \infty} \left[ C ( \textnormal{e}^{x(t-v)} + 1 ) - 2 C (\textnormal{e}^{x(t-v)} - 1) \right.  \\
 & \quad \; \left. - C\frac{|y|}{x} (\textnormal{e}^{x(t-v)} - 1) \right] \\
 & = C \left(3 + \frac{|y|}{x} \right).
\end{align*}
\end{proof}

\begin{lemma}\label{lemma3}
Let $X$ be a Banach space, $A: \mathbb{R_+} \to X$ locally of bounded variation, $A(0)=0$ and suppose that

\begin{equation*}
 \sup\limits_{t \in \mathbb{R_+}} \left| \left| \textnormal{e}^{-x_0t} \int_0^{t} \textnormal{e}^{x_0s} \textnormal{d}A(s) \right| \right| \leq C
\end{equation*}
for some $x_0>0$ and some $C>0$. Then 

\begin{equation*}\label{l3}
 \sup\limits_{t \in \mathbb{R_+}} \left| \left| \textnormal{e}^{-xt} \int_0^{t} \textnormal{e}^{xs} \textnormal{d}A(s) \right| \right| \leq \frac{Cx_0}{x}
\end{equation*}
for every $x$ with $0<x \leq x_0$.
\end{lemma}
\vspace{0.1cm}
\begin{proof}
By using $G(s) = \int_0^s \textnormal{e}^{x_0r} \textnormal{d}A(r), s \in \bR_+$ we have

\begin{align*}
\left| \left| \int_0^t \textnormal{e}^{xs} \textnormal{d}A(s) \right| \right| 
 & = \left| \left| \int_0^t \textnormal{e}^{xs} \textnormal{e}^{-x_0s} \te^{x_0s} \textnormal{d}A(s) \right| \right| \\
 & = \left| \left| \int_0^t \textnormal{e}^{s(x-x_0)} \textnormal{d} G(s) \right| \right|.
\end{align*}
Remember that $||G(s)|| \leq C\te^{x_0s}$ according to the assumptions. We integrate by parts and estimate:

\begin{align}
 \left| \left| \int_0^t \textnormal{e}^{s(x-x_0)} \textnormal{d} G(s) \right| \right| 
 & = \left| \left| \left[ \textnormal{e}^{s(x-x_0)} G(s) \right]_0^t - (x-x_0) \int_0^t \textnormal{e}^{s(x-x_0)} G(s) \textnormal{d}s \right| \right|  \nonumber \\ 
 & \leq \left| \left| \textnormal{e}^{t(x-x_0)} G(t) \right| \right| + | x-x_0 | \int_0^t \te^{xs} \te^{-x_0s}  || G(s) || \textnormal{d}s \nonumber \\ 
 & \leq C \textnormal{e}^{xt} + C | x-x_0 | \int_0^t \textnormal{e}^{xs} \textnormal{d}s \nonumber \\
 & = C \textnormal{e}^{xt} + C \frac{| x-x_0 |}{x} ( \textnormal{e}^{xt} - 1) \nonumber \label{lst} \\ 
 & \leq \left( 1 + \frac{| x-x_0|}{x}  \right) C \textnormal{e}^{xt}. 
\end{align}
The coefficient can be simplified to

\[
 1 + \frac{x_0 - x}{x} = \frac{x_0}{x} 
\]
and we get the result.
\end{proof}

\begin{bem}\label{bem}
\rm For proving Theorem \ref{TH} it would be nice to extend the result of Lemma \ref{lemma3} to all $x \in (0,R(t)]$. In fact, this is not possible for all 
$t \in \mathbb{R}_+$. For example, fix $T>0$ and consider $A: \bR_+ \to \{0,1\}$ with

\[
 A(t) = 
\left\{
\begin{array}{ll}
 0, & 0 \leq t \leq T, \\
 1, & t > T.
\end{array}
\right.
\] 
Then $A$ is of bounded variation and $A(0)=0$. Let $x>0$. Define $g_{x}: \bR_+ \to \bR_+$ by

\[
g_{x}(t) \coloneqq \te^{-xt} \int_0^t \te^{xt} \td A(t) = 
\left\{
\begin{array}{ll}
 0, & 0 \leq t \leq T, \\
 \te^{x(T-t)}, & t>T.
\end{array}
\right.
\]
Now choose $x_0=1$ and set $C \coloneqq 1 = \sup\limits_{t \in \bR_+} |g_{x_0}(t)|$. Because $\sup\limits_{t \in \bR_+} |g_{x}(t)| = 1$ holds for every 
$x > 0$ we have

\[
\sup\limits_{t \in \bR_+} |g_{x}(t)| > \frac{1}{x} = \frac{Cx_0}{x}
\]
for every $x>x_0$. But notice that there exists $T'>T$ so that 

\[
 \sup\limits_{t>T'} |g_{x}(t)| < \frac{1}{x} = \frac{Cx_0}{x} 
\] 
is true for all $x>0$ and, in particular, for all $x \in (0,R(t)]$.

\end{bem}

\section{Proof of Theorem \ref{TH}}

First, we show that the condition 

\begin{equation}\label{tc}
 \sup\limits_{t \in \mathbb{R_+}} \left| \left| \textnormal{e}^{-x_0t} \int_0^{t} \textnormal{e}^{x_0s} \textnormal{d}A(s) \right| \right| \leq C
\end{equation}
for some $x_0>0$ and some $C>0$ is sufficient for the existence of the Laplace-Stieltjes transform $f(z) = \int_0^{\infty} \textnormal{e}^{-zs} \textnormal{d}A(s)$
of $A$ for every $z \in H$. Using inequality \eqref{lst} in the proof of Lemma \ref{lemma3} we conclude that \eqref{tc} is vaild for every $x>0$. By Lemma \ref{lemma2}
we get that $\int_t^{\infty} \textnormal{e}^{-zs} \textnormal{d}A(s)$ exists for every $z = x + \ti y$ with $x>0$ and $y \in \bR$. Therefore, we proved the above claim.  

For proving the quantitative statement we use the notation

\[
 f_t(z) = \int_0^t \textnormal{e}^{-zs}  \textnormal{d}A(s)
\]
so that

\[
 f_t(0) = \int_0^t \textnormal{d}A(s) = A(t)
\]
and consider now the behaviour of $||f_t(0) - f(0)||$ for large $t$. Fix $t>T$ and consider $R \in [1,R(t)]$. We form a contour $\Gamma$ which consists of two parts:
$\Gamma_1$ is the arc $\{z \in \mathbb{C} \, | \, |z|=R, \textup{Re}(z) \geq 0 \}$ in the closed right half-plane. $\Gamma_2$ consists of the three segments 
$[\ti R,-\frac{1}{2M(R)} + \ti R]$, $[-\frac{1}{2M(R)} + \textnormal{i}R,- \frac{1}{2M(R)} - \ti R]$ and $[- \frac{1}{2M(R)} - \ti R, -\ti R]$. Therefore, $\Gamma_2$
is contained in $Q$. By Cauchy's integral formula we get

\begin{equation}\label{int} 
||f_t(0) - f(0)|| = \left| \left| \frac{1}{2\pi \textnormal{i}} \int_{\Gamma} \frac{f_t(z)-f(z)}{z} \textnormal{e}^{tz} \left( 1 + \frac{z^2}{R^2} \right)^2 \textnormal{d}z \right| \right|.
\end{equation}
As $f_t$ is an entire function we can replace the integral $\int_{\Gamma_2} f_t(z) \frac{\textnormal{e}^{tz}}{z} ( 1 + \frac{z^2}{R^2} )^2 \textnormal{d}z$ 
by an integral over $\tilde{\Gamma}_1 = \{z \in \mathbb{C} \, | \, |z|= R, \textup{Re}(z) < 0 \}$, which is the reflection of $\Gamma_1$ through the origin. 
Let us split the integral in \eqref{int} into three parts:

\begin{align*}
||f_t(0) - f(0)|| & \leq \left| \left| \frac{1}{2\pi \textnormal{i}} \int_{\Gamma_1} \frac{f_t(z)-f(z)}{z} \textnormal{e}^{tz} \left( 1 + \frac{z^2}{R^2} \right)^2 \textnormal{d}z \right| \right| \\
                  & \quad \;  + \left| \left| \frac{1}{2\pi \textnormal{i}} \int_{\tilde{\Gamma}_1} \frac{f_t(z)}{z} \textnormal{e}^{tz} \left( 1 + \frac{z^2}{R^2} \right)^2 \textnormal{d}z \right| \right| \\
                  & \quad \;  + \left| \left| \frac{1}{2\pi \textnormal{i}} \int_{\Gamma_2} \frac{f(z)}{z} \textnormal{e}^{tz} \left( 1 + \frac{z^2}{R^2} \right)^2 \textnormal{d}z \right| \right| \\  
                  & \eqqcolon I + II + III
\end{align*}
Now, we estimate every single integral. For that we use 

\[
 \left| 1 + \frac{z^2}{R^2} \right| = \frac{2|x|}{R}, \quad \frac{1}{|z|} = \frac{1}{R}, \quad |y| \leq R
\]
on the circle $|z|=R$, with $z = x + \ti y$. By Lemma \ref{lemma3} we know that 

\begin{equation*}
\sup\limits_{t > T} \left|\left| \textnormal{e}^{-xt} \int_0^{t} \textnormal{e}^{xs} \textnormal{d}A(s) \right|\right| \leq \frac{C}{x_0} \cdot \frac{x_0}{x} = \frac{C}{x}
\end{equation*}
for every $x$ with $0<x<x_0$. So all estimations do not depend on whether $x$ is smaller than $x_0$ or fulfills $x_0 \leq x \leq R$. 
Then $I$ can be estimated by (remember $z=x+\textnormal{i}y$ and notice that $x \geq 0$ on $\Gamma_1$)

\begin{align*}
 I & \leq \frac{1}{2\pi} \int_{\Gamma_1} ||(f_t(z)-f(z))|| \textnormal{e}^{tx}  \frac{4x^2}{R^3} |\textnormal{d}z| \\
   & = \frac{1}{2\pi} \int_{\Gamma_1} \left| \left| \textnormal{e}^{tx} \int_t^{\infty} \textnormal{e}^{-zs} \textnormal{d}A(s) \right| \right| \frac{4x^2}{R^3} |\textnormal{d}z| \\
   & \leq \frac{1}{2\pi} \int_{\Gamma_1} \frac{C}{x} \left( 3 + \frac{|y|}{x} \right) \frac{4x^2}{R^3} |\textnormal{d}z| \\
   & \leq \frac{6C}{\pi R^3} \int_{\Gamma_1} x \; |\textnormal{d}z| + \frac{2C}{\pi R^2} \int_{\Gamma_1} |\textnormal{d}z| \\
   & = \frac{12C}{\pi R} + \frac{2C}{R} \leq \frac{6C}{R}, 
\end{align*}
where we used Lemma \ref{lemma2}. For the estimation of $II$ we assume $x < 0$ and define $\tilde{z} \coloneqq -z = -x-\ti y$, with $y \in \bR$ and $\tilde{z}$ 
lies in the right half-plane. Now, similar considerations give  

\begin{align*}
II & \leq \frac{1}{2\pi} \int_{\tilde{\Gamma}_1} ||f_t(z)|| \textnormal{e}^{tx}  \frac{4(-x)^2}{R^3} |\textnormal{d}z| \\
   & = \frac{1}{2\pi} \int_{\tilde{\Gamma}_1} \left| \left| \textnormal{e}^{tx} \int_0^t \textnormal{e}^{\tilde{z}s} \textnormal{d}A(s) \right| \right| \frac{4(-x)^2}{R^3} |\textnormal{d}z| \\
   & \leq \frac{1}{2\pi} \int_{\tilde{\Gamma}_1} \frac{C}{-x} \left( 1 + \frac{|y|}{-x} \right) \frac{4(-x)^2}{R^3} |\textnormal{d}z| \\
   & \leq \frac{2C}{\pi R^3} \int_{\tilde{\Gamma}_1} (-x) \; |\textnormal{d}z| + \frac{2C}{\pi R^2} \int_{\tilde{\Gamma}_1} |\textnormal{d}z| \\
   & = \frac{4C}{\pi R} + \frac{2C}{R} \leq \frac{4C}{R},
\end{align*}
where we used Lemma \ref{lemma1}. Finally, we consider $III$. By assumption, $||f(z)||$ is less than or equal to $M(|y|)$ for every $z=x+\textnormal{i}y$ on the path of 
integration. Along the segments $[\ti R,-\frac{1}{2M(R)} + \ti R]$ and $[- \frac{1}{2M(R)} - \ti R, -\ti R]$ we have (remember $R \geq 1, M(R) \geq 1$)

\[
 \left| 1 + \frac{z^2}{R^2} \right| \leq \frac{\sqrt{2}}{R} \; \; \textnormal{and} \; \; \frac{1}{|z|} \leq \frac{1}{R}.
\]
For the segment from $-\frac{1}{2M(R)} + \ti R$ to $-\frac{1}{2M(R)} - \ti R$ we can estimate

\[
 \left| 1 + \frac{z^2}{R^2} \right|  \leq \sqrt{2} \; \; \textnormal{and} \; \; \frac{1}{|z|} \leq 2 M(R).
\]
Therefore, we get

\begin{align*}
 III & \leq \frac{1}{2\pi} \left( 2 \int_{-\frac{1}{2M(R)}}^0 M(R) \textnormal{e}^{tx} \frac{1}{R} \cdot \frac{2}{R^2} \textnormal{d}x \right) \\
     & \quad \; + \frac{1}{2\pi} \int_{-R}^R M(R) \textnormal{e}^{-\frac{t}{2M(R)}} \cdot 2M(R) \cdot 2 \; \textnormal{d}y. 
\end{align*}
Since 

\[
 \int_{-\frac{1}{2M(R)}}^0 e^{tx} \textnormal{d}x = \frac{1}{t} - \frac{1}{t} \textnormal{e}^{-\frac{t}{2M(R)}} \leq \frac{1}{t},  
\] 
we conclude

\begin{equation*}
 III \leq \frac{M(R)}{t R^3} + 2 R (M(R))^2 \textnormal{e}^{-\frac{t}{2M(R)}}. 
\end{equation*}
If we summarize all estimates we have

\begin{equation}\label{Abs}
 || A(t) - f(0) ||  \leq \frac{10C}{R} + \frac{M(R)}{t R^3} + 2 R (M(R))^2 \textnormal{e}^{-\frac{t}{2M(R)}}.
\end{equation}  
Now, we optimize this estimate over $R$ by equating the first and the third term of the right-hand side:

\begin{equation*}
 \frac{10C}{R_{\opt}} = 2 R_{\opt} (M(R_{\opt}))^2 \textnormal{e}^{-\frac{t}{2M(R_{\opt})}}.
\end{equation*}  
Therefore, we get 

\begin{equation}\label{ropt}
 t = 4 M(R_\opt)\left( \log R_\opt + \log M(R_\opt) - \frac{1}{2} \log(5C) \right),    
\end{equation}
that is

\begin{equation*}
 R_{\opt} = M_{\log}^{-1} \left( \frac{t}{4} \right),
\end{equation*}
where $M^{-1}_{\log}$ is the inverse of the function on the right-hand side of \eqref{ropt}, i.e. 
$M_{\log}(\cdot) = M(\cdot)( \log \cdot + \log M(\cdot) - \frac{1}{2} \log(5C))$. 

Since $R \geq 1$ by assumption, we have $t \geq 4 M(1)\left( \log M(1) - \frac{1}{2} \log (5C) \right)$. So we define 

\[
T' \coloneqq \max \; \left\{ T \; , \; 4 M(1) \left( \log M(1) - \frac{1}{2} \log (5C) \right) \right\}. 
\]
If we insert $t$ into the middle term of the sum in \eqref{Abs} it follows

\[
 \frac{M(R_\opt)}{t R_\opt^3} \leq \frac{1}{\log M(1) - \log ( \sqrt{5C} )} \cdot \frac{1}{R_\opt^3} \leq \frac{K'}{R_\opt} = \frac{K'}{M_{\log}^{-1} \left( \frac{t}{4} \right)} , 
\]
where $K' \coloneqq (\log M(1) - \log ( \sqrt{5C} ))^{-1}$. Finally, we check if $R_\opt \in [1,R(t)]$ is true. If $R_\opt \in [1,R(t)]$, then
\begin{equation}\label{1est}
  || A(t) - f(0) ||  \leq \frac{20C}{R_\opt} + \frac{K'}{R_\opt} \leq \frac{K}{M^{-1}_{\log} \left(\frac{t}{4} \right)}
\end{equation}
for every $t > T'$, with a suitable $K>0$. Otherwise, we have

\begin{equation*}
 \frac{1}{M_{\log}^{-1} \left( \frac{t}{4} \right)} = \frac{1}{R_\opt} \leq \frac{1}{R(t)},
\end{equation*}
so that

\begin{equation}\label{2est}
  || A(t) - f(0) ||  \leq \frac{20C}{R_\opt} + \frac{K'}{R_\opt} \leq \frac{K}{R(t)}
\end{equation}
for every $t > T'$, with a suitable $K>0$. Combining \eqref{1est} and \eqref{2est} gives the result of Theorem \ref{TH}.

\begin{bem}
\rm The integrand in \eqref{int} is multiplied by the terms $\left( 1 + \frac{z^2}{R^2} \right)^2$, which is the so-called "fudge factor", and $\te^{tz}$. Both do not 
change the value of the contour integral but help to estimate the integral, see \cite{Korevaar82}. This idea is due to Newman \cite{Newman80}. 
\end{bem}

\section{Different Remarks}\label{remarks}
In this section we show that Theorem \ref{TH} includes the results from \cite{Batty08} and can be applied to Dirichlet series with bounded coefficients.

\begin{bem}\label{bounded functions}
\rm Let $X$ be a Banach space and $f: \bR_+ \to X$ a bounded measurable function. We know that $\hat{f}(z) = \int_0^{\infty} \te^{zs} f(s) \td s$
exists for every $z \in H$ and that $f \in L^1_{loc}(\bR_+;X)$. Furthermore, the function $A: \bR_+ \to X$ with $A(t) \coloneqq \int_0^t f(s) \td s$ is locally of 
bounded variation and differentiable almost everywhere. So $A$ is an antiderivative of $f$ with $A(0)=0$. Denote by $C$ the bound of $||f(s)||$.
Then 

\begin{align*}
 \left|\left| \textnormal{e}^{-xt} \int_0^{t} \textnormal{e}^{xs} \textnormal{d}A(s) \right|\right|
 & = \left| \left| \te^{-xt} \int_0^{t} \te^{xs} f(s) \td s \right| \right| \\
 & \leq \te^{-xt} \int_0^{t} \te^{xs} || f(s) || \td s \\
 & \leq C \te^{-xt} \cdot \frac{1}{x} (\te^{xt} - 1) \\
 & \leq \frac{C}{x},
\end{align*}
for every $t \in \bR_+$ and every $x>0$. Choosing $R(t) = \infty$ for every $t>0$ we see that condition \eqref{TC} is satisfied. For this choice of $R(t)$, the
condition $R_\opt \in [1,R(t)]$ is always true. 

Further, if we make the same assumptions as in Section 4 of \cite{Batty08}, all conditions of our Theorem \ref{TH} are fulfilled. So we can conclude that
there exist $T' \geq 0, K>0$ so that

\begin{equation}\label{estimate}
 \left| \left| \int_0^t f(s) \td s - f(0) \right| \right|  \leq K \max \left\{ \frac{1}{M^{-1}_{\log} \left( \frac{t}{4} \right)}, \frac{1}{R(t)} \right\}
\end{equation}
for every $t>T'$, where $M^{-1}_{\log}$ is the inverse of the function $M_{\log}$ defined by $M_{\log}(a) = M(a)( \log a + \log M(a) - \frac{1}{2} \log(5C))$ 
for $a \geq 1$. 

Note that by the above choice of $R(t)$ the maximum in \eqref{estimate} is equal to $(M^{-1}_{\log}(\frac{t}{4}))^{-1}$ for every $t>T'$. 
In this way we recover \cite[Theorem 4.1.]{Batty08} in the case $k=1$.
\end{bem}

\begin{bem}
\rm We show that condition \eqref{TC} of Theorem \ref{TH} is automatically true for special Dirichlet series.  

Let $X$ be a Banach space and $(b_n)_{n \in \mathbb{N}} \in l^{\infty}(X)$. Set $D \coloneqq \max\{||b||_{\infty},1\}$. Define a sequence 
$(a_n)_{n \in \mathbb{N}}$ by $a_n \coloneqq \frac{b_n}{n}$ for every $n \in \mathbb{N}$. Consider the Dirichlet series

\[
 f(z) = \sum_{n=1}^{\infty} \frac{a_n}{n^z} = \sum_{n=1}^{\infty} \frac{b_n}{n^{z+1}}, 
\]
which is analytic in $H$. Furthermore, we define $A: \bR_+ \to X$ by $A(s) \coloneqq \sum_{\log n < s} a_n$ so that $A$ is
locally of bounded variation, continuous from the left and $A(0)=0$. 
For $t>0$ and $x>0$ we get

\begin{align*}
 \left|\left| \textnormal{e}^{-xt} \int_0^{t} \textnormal{e}^{xs} \textnormal{d}A(s) \right|\right|
 & = \left| \left| \te^{-xt} \sum_{\log n \, < \, t} \te^{x\log n} a_n \right| \right| \\
 & = \left| \left| \te^{-xt} \sum_{\log n \, < \, t} n^x a_n \right| \right| \\
 & \leq D \te^{-xt} \sum_{\log n \, < \, t} n^{x-1} \\
 & \leq D \te^{-xt} \sum_{n=1}^{\floor{\te^t}} n^{x-1}. 
\end{align*}
Further estimates yield

\begin{align*}
 \te^{-xt} \sum_{n=1}^{\floor{\te^t}} n^{x-1} & \leq \te^{-xt} \int_0^{\te^t} (s+1)^{x-1} \td s \\
                                              & \leq \frac{1}{x} \te^{-xt} (e^t + 1)^x \\
					      & = \frac{1}{x} (1 + \te^{-t})^x, 
\end{align*}
which is bounded by $\te \cdot x^{-1}$ for every $t \in \bR_+$ and every $x$ with $0 < x \leq \te^t$. Define $R(t) \coloneqq \te^t$  
for every $t>0$. It follows that $R(t)>1$. Defining $C \coloneqq D \cdot \te$ we conclude

\[
 \left|\left| \textnormal{e}^{-xt} \int_0^{t} \textnormal{e}^{xs} \textnormal{d}A(s) \right|\right| \leq \frac{C}{x}
\]
for every $t>0$ and $x \in (0,\te^t]$. Therefore, we state the following corollary of Theorem \ref{TH}.
\end{bem}

\vspace{0.3cm}
\begin{kol}
Let $X$ be a Banach space and $(b_n)_{n \in \mathbb{N}} \in l^{\infty}(X)$. Define $a_n \coloneqq \frac{b_n}{n}$ for every $n \in \mathbb{N}$. Then the 
Dirichlet series

\[
 f(z) = \sum_{n=1}^{\infty} \frac{a_n}{n^{z}}
\]
is analytic in $\{ z \in \bC \; | \; \textup{Re}(z) > 0\}$. Let $M: [0,\infty) \to [1,\infty)$ be a continuous, increasing function and define 
$R: [0,\infty) \to [1,\infty]$ by $R(t) \coloneqq \te^t$. Assume that $f$ has an analytic extension into the region

\begin{equation*}
  Q \coloneqq \left\{ z \in \mathbb{C} \; \big| \; 0 \geq x > - \frac{1}{M(|y|)} \right\}
\end{equation*}
so that $||f(z)|| \leq M(|y|) $ holds throughout $Q$, where $z=x+\ti y$. Then there exist $K>0, T'\geq 0$ so that

\[
 \bigg|\bigg| \sum_{\log n < t} a_n - f(0) \bigg|\bigg|  \leq K \max \left\{ \frac{1}{M^{-1}_{\log} \left( \frac{t}{4} \right)}, \frac{1}{R(t)} \right\}
\]
for every $t>T'$, where $M^{-1}_{\log}$ is the inverse of the function $M_{\log}$ defined by $M_{\log}(a) = M(a)( \log a + \log M(a) - \frac{1}{2} \log(5C))$ 
for $a \geq 1$. 
\end{kol}

\vspace{0.3cm}
\subsection*{Acknowledgement}
I am grateful to my supervisor Ralph Chill for his support and his suggestions, which helped to improve the paper a lot. 

\bibliographystyle{plain}

\end{document}